\documentclass[12]{amsproc}

\usepackage{amssymb}

\usepackage[dvipdfm,bookmarks=true,bookmarksnumbered=true,bookmarkstype=toc]{hyperref}   
\usepackage{amsmath, amsthm, amsfonts}
\usepackage{amscd}
\usepackage{verbatim}

\numberwithin{equation}{section}

\numberwithin{figure}{section}

\newtheorem{thm}{Theorem}\numberwithin{thm}{section}
\newtheorem{lem}{Lemma}\numberwithin{lem}{section}
\newtheorem{corollary}{Corollary}\numberwithin{corollary}{section}
\newtheorem{conjecture}{Conjecture}\numberwithin{conjecture}{section}
\newtheorem{question}{Question}\numberwithin{question}{section}
\numberwithin{assertion}{section}
\theoremstyle{definition}
\newtheorem{defn}{Definition}\numberwithin{defn}{section}
\theoremstyle{remark}
\newtheorem{rem}{Remark}\numberwithin{rem}{section}
\newtheorem{example}{Example}\numberwithin{example}{section}
\newcommand{\be}{\begin{equation}}
\newcommand{\beq}{\begin{equation}}
\newcommand{\eeq}{\end{equation}}
\newcommand{\beqa}{\begin{eqnarray}}
\newcommand{\eeqa}{\end{eqnarray}}
\newcommand{\beqanl}{\begin{eqnarray*}}
\newcommand{\eeqanl}{\end{eqnarray*}}
\newcommand{\bs}{\begin{sub}}
\newcommand{\es}{\end{sub}}
\newcommand{\bsn}{\begin{subn}}
\newcommand{\esn}{\end{subn}}
\newcommand{\bea}{\begin{eqnarray}}
\newcommand{\eea}{\end{eqnarray}}
\newcommand{\bean}{\begin{eqnarray*}}
\newcommand{\eean}{\end{eqnarray*}}
\newcommand{\BA}[1]{\begin{array}{#1}}
\newcommand{\EA}{\end{array}}

\newlength{\wex}  \newlength{\hex}
            \def\gl{\lambda}

\def\squarebox#1{\hbox to #1{\hfill\vbox to #1{\vfill}}}

\newcommand{\Green}[4]{\mbox{$G^{#1}_{#2}(#3,#4)$}}

\title[Davies' conjecture]
{On Davies' conjecture and strong ratio limit
properties for the heat kernel}

\author[Y.~Pinchover]
{Yehuda Pinchover}

\address{Department of Mathematics, Technion - Israel Institute of
Technology,  Haifa, Israel} \email{pincho@techunix.technion.ac.il}

\thanks{Partially supported by the Israel Science Foundation founded by the
Israeli Academy of Sciences and Humanities.}

\begin{document}

\begin{abstract}
We study strong ratio limit properties and the exact long time
asymptotics of the heat kernel of a general second-order parabolic
operator which is defined on a noncompact Riemannian manifold.
\end{abstract}

\maketitle

\section{Introduction}\label{Introduction}

 Let $P$ be a linear, second-order, elliptic operator defined on a  noncompact, connected,
$C^3$-smooth Riemannian manifold $\mathcal{M}$ of dimension $d$.
Here $P$ is an  elliptic operator with real, H\"{o}lder continuous
coefficients which in any coordinate system
$(U;x_{1},\ldots,x_{d})$ has the form
$$P(x,\partial_{x})=-\sum_{i,j=1}^{d}
a_{ij}(x)\partial_{i}\partial_{j} + \sum_{i=1}^{d}
b_{i}(x)\partial_{i}+c(x).$$  \noindent
 We assume that
for each $x\in \mathcal{M}$ the real quadratic form
$\sum_{i,j=1}^{d}\! a_{ij}(x)\xi_{i}\xi_{j}$ is positive definite.
The formal adjoint of $P$ is denoted by $P^*$. Denote the cone of
all positive (classical) solutions of the equation $Pu\!=\!0$ in
$\mathcal{M}$ by $\mathcal{C}_{P}(\mathcal{M})$. The {\em
generalized principal eigenvalue}  is defined by
$$\gl_0=\gl_0(P,\mathcal{M})
:= \sup\{\gl \in \mathbb{R} \; :\;
\mathcal{C}_{P-\lambda}(\mathcal{M})\neq \emptyset\}.$$ \noindent
 Throughout this paper we always assume that
$\lambda_0\geq 0$ (actually, as it will become clear below, it is
enough to assume that $\lambda_0>-\infty$).

We consider the parabolic operator $L$
\begin{equation}\label{eqL}
  Lu=u_t+Pu \qquad \mbox{ on } \mathcal{M}\times (0,\infty).
\end{equation}  \noindent
 We denote by $\mathcal{H}_P(\mathcal{M}\times (a,b))$ the cone of all
nonnegative solutions of the equation $Lu=0$ in $\mathcal{M}\times
(a,b)$. Let $k_P^{\mathcal{M}}(x,y,t)$ be the {\em minimal
(positive) heat kernel} of the parabolic operator $L$ in
$\mathcal{M}$. If for some $x\neq y$
$$\int_0^\infty k_P^{\mathcal{M}}(x,y,\tau)\,\mathrm{d}\tau<\infty \qquad
\left(\mbox{respect., } \int_0^\infty
k_P^{\mathcal{M}}(x,y,\tau)\,\mathrm{d}\tau=\infty\right),$$ then
$P$ is said to be a {\em subcritical} (respect., {\em critical})
operator in $\mathcal{M}$ \cite{Pinsky}.

Recall that if $\lambda\!<\!\lambda_0$, then $P\!-\!\lambda$ is
subcritical in $\mathcal{M}$, and for $\lambda\leq \lambda_0$, we
have $k_{P\!-\!\lambda}^\mathcal{M}(x,y,t)\!=\!\mathrm{e}^{\lambda
t}k_P^\mathcal{M}(x,y,t)$. Furthermore, $P$ is critical (respect.,
subcritical) in $\mathcal{M}$, if and only if $P^*$ is critical
(respect., subcritical) in $\mathcal{M}$. If $P$ is critical in
$\mathcal{M}$, then there exists a unique positive solution
$\varphi \in\mathcal{C}_{P}(\mathcal{M})$ satisfying
$\varphi(x_0)=1$. This solution is called the {\em ground state of
the operator $P$ in $\mathcal{M}$} \cite{Pheat,Pinsky}. The ground
state of $P^*$ is denoted by $\varphi^*$. A critical operator $P$
is said to be {\em positive-critical} in $\mathcal{M}$ if
$\varphi^*\varphi\in L^1(\mathcal{M})$, and {\em null-critical} in
$\mathcal{M}$ if $\varphi^*\varphi\not\in L^1(\mathcal{M})$. In
\cite{Pheat,P03} we proved:
\begin{thm}\label{mainthm}
Let $x,y\in \mathcal{M}$. Then
 \begin{equation*}
\lim_{t\to\infty} \mathrm{e}^{\lambda_0
t}k_P^{\mathcal{M}}(x,y,t)\!=\!
  \begin{cases}
    \dfrac{\varphi(x)\varphi^*(y)}{\int_\mathcal{M}\!
\varphi(z)\varphi^*(z)\,\mathrm{d}z} & \text{if } P\!-\!\lambda_0
\text{ is positive-critical},
\\[3mm]
    0 & \text{otherwise}.
  \end{cases}
 \end{equation*}
Furthermore, for $\lambda<\lambda_0$, let
$\Green{\mathcal{M}}{P-\lambda}{x}{y}\!:=\!\int_0^\infty \!
k_{P-\lambda}^{\mathcal{M}}(x,y,\tau)\mathrm{d}\tau$ be the
minimal (positive) Green function of the operator $P\!-\!\lambda$
on $\mathcal{M}$. Then
\begin{equation}\label{eqgreen}
\lim_{t\to\infty} \mathrm{e}^{\lambda_0
t}k_P^{\mathcal{M}}(x,y,t)=
\lim_{\lambda\nearrow\lambda_0}(\lambda_0-\lambda)\Green{\mathcal{M}}{P-\lambda}{x}{y}.
\end{equation}
\end{thm}
Having proved that $\lim_{t\to\infty} \mathrm{e}^{\lambda_0
t}k_P^{\mathbf{M}}(x,y,t)$ always exists, we next ask how fast
this limit is approached.  It is natural to conjecture that the
limit is approached equally fast for different points $x,y\in
\mathcal{M}$. Note that in the context of Markov chains, such an
{\em (individual) strong ratio limit property} is in general not
true \cite{Chu}. The following conjecture was raised by
E.~B.~Davies \cite{D} in the selfadjoint case.
\begin{conjecture}\label{conjD}
Let $Lu=u_t+P(x, \partial_x)u$ be a parabolic operator which is
defined on a Riemannian manifold $\mathcal{M}$. Fix a reference
point $x_0\in \mathcal{M}$. Then
\begin{equation}\label{eqconjD}
\lim_{t\to\infty}\frac{k_P^\mathcal{M}(x,y,t)}{k_P^\mathcal{M}(x_0,x_0,t)}=a(x,y)
\end{equation}
exists and is positive for all $x,y\in \mathcal{M}$.
\end{conjecture}
\begin{rem}\label{6Rem82}
Theorem~\ref{mainthm} implies that Conjecture \ref{conjD} holds
true in the positive critical case. So, we may assume in the
sequel that $P$ {\bf is not positive critical}. Also, Conjecture
\ref{conjD} does not depend on the value of $\lambda_0$, hence
from now on, {\bf we shall assume that $\lambda_0=0$}.
\end{rem}
\begin{rem}\label{6Rem29}
In the selfadjoint case, Conjecture \ref{conjD} holds true  if
$\dim \mathcal{C}_{P}(\mathcal{M}) = 1$ \cite[Corollary 2.7]{ABJ}.
In particular, it holds true for a critical selfadjoint operator.
Therefore, it would be interesting to prove Conjecture \ref{conjD}
at least under the assumption
\begin{equation}\label{eqoned}
\dim \mathcal{C}_{P}(\mathcal{M})=\dim
\mathcal{C}_{P^*}(\mathcal{M})=1,
\end{equation}
which holds true in the critical case and in many important
subcritical cases. Recently,  Agmon \cite{Agp} has obtained the
exact asymptotics (in $(x,y,t)$) of the heat kernel for a periodic
(non-selfadjoint) operator on $\mathbb{R}^d$. It follows from
Agmon's results that Conjecture \ref{conjD} holds true in this
case. For a probabilistic interpretation of Conjecture
\ref{conjD}, see \cite{ABJ}.
\end{rem}
\begin{rem}\label{6Rem21}
Let $t_n\to \infty$. By a standard parabolic argument, we may
extract a subsequence $\{t_{n_k}\}$ such that for every $x,y\in
\mathcal{M}$ and $s<0$
$$a(x,y,s):=\lim_{k\to\infty}\frac{k_P^\mathcal{M}(x,y,s+t_{n_k})}{k_P^\mathcal{M}(x_0,y_0,t_{n_k})}
$$ exists. Moreover, $a(\cdot,y,\cdot)\in
\mathcal{H}_P(\mathcal{M}\times \mathbb{R}_-)$. Note that in the
selfadjoint case, the above is valid for all $s\in \mathbb{R}$.
 \end{rem}
 \begin{rem}\label{6Rem2}
The example constructed in \cite[Section 4]{PSt} shows a case
where Conjecture \ref{conjD} holds true on $\mathcal{M}$, while
the limit function $a(x,y)=1$  is not a $\lambda_ 0$-invariant
positive solution. Compare this with \cite[Theorem 25]{D} and the
discussion therein above Lemma 26. Note also that in general, the
limit function $a(x,y)$ in (\ref{eqconjD}) need not be a {\em
product} of solutions of the equations $Pu=0$ and $P^*u=0$, as is
demonstrated in \cite{CMS1}, in the hyperbolic space, and in
Example \ref{examinim1}.
 \end{rem}

\section{Strong ratio properties}
In the symmetric case the function $t\mapsto
k_P^\mathcal{M}(x,x,t)$ is log-convex, and therefore, a
polarization argument implies  that $\lim_{t\to\infty}
\frac{k_P^\mathcal{M}(x,y,t+s)}{k_P^\mathcal{M}(x,y,t)}=1$ for all
$x,y\in \mathcal{M}$ and $s\in \mathbb{R}$ \cite{D}. In the
nonsymmetric case we have.
\begin{lem}\label{assliminfsup}
For every  $x,y\in \mathcal{M}$ and $s\in \mathbb{R}$, we have
that
\begin{equation}\label{eqskeleton11}
\liminf_{t\to\infty}
\frac{k_P^\mathcal{M}(x,y,t+s)}{k_P^\mathcal{M}(x,y,t)}\leq 1\leq
\limsup_{t\to\infty}
\frac{k_P^\mathcal{M}(x,y,t+s)}{k_P^\mathcal{M}(x,y,t)}\,.
\end{equation}
Similarly, for any $s>0$
\begin{equation}\label{eqskeleton19}
\liminf_{n\to\infty}
\frac{k_P^\mathcal{M}(x,y,(n\!\pm\!1)s)}{k_P^\mathcal{M}(x,y,ns)}\leq
1\leq\limsup_{n\to\infty}
\frac{k_P^\mathcal{M}(x,y,(n\!\pm\!1)s)}{k_P^\mathcal{M}(x,y,ns)}\,.
\end{equation}
In particular, if $\lim_{t\to\infty}
[k_P^\mathcal{M}(x,y,t\!+\!s)/k_P^\mathcal{M}(x,y,t)]$ exists, it
equals to $1$.
\end{lem}
\begin{proof}
 We may assume that $P$ is not positive-critical. Let $s<0$. By
Theorem~\ref{mainthm} and the parabolic Harnack inequality we have
\begin{equation}\label{eqskeleton1}
1\leq \limsup_{t\to\infty}
\frac{k_P^\mathcal{M}(x,y,t+s)}{k_P^\mathcal{M}(x,y,t)}\leq C(s,y)
.
\end{equation}
Suppose that $\liminf_{t\to\infty}
\frac{k_P^\mathcal{M}(x,y,t+s)}{k_P^\mathcal{M}(x,y,t)}=\ell > 1$.
It follows that there exists $0<q<1$ and $T_s>0$ such that
$$k_P^\mathcal{M}(x,y,t)< q k_P^\mathcal{M}(x,y,t+s)\qquad \forall T_s<t.$$
 By induction and the Harnack inequality, we obtain
that there exist $\mu<0$ and $C>0$ such that
$k_P^\mathcal{M}(x,y,t)<C\mathrm{e}^{\mu t}$ for all $t>1$, a
contradiction to  the assumption $\lambda_0\!=\!0$. Therefore,
(\ref{eqskeleton11}) is proved for $s<0$, which in turn implies
(\ref{eqskeleton11}) also for $s>0$.  (\ref{eqskeleton19}) can be
proven similarly.
 \end{proof}
\begin{rem} The condition $\liminf
_{t\to\infty}\frac{k_P^\mathcal{M}(x,y,t+s)}{k_P^\mathcal{M}(x,y,t)}\geq
1$ for $s>0$ is sometimes called {\em Lin's  condition}
\cite{LM}.\end{rem}
\begin{corollary}\label{corskeleton}
Let $x, y \in \mathcal{M}$. Suppose that
\begin{equation}\label{eqliminf5}
\lim_{n\to\infty}
\frac{k_P^\mathcal{M}(x,y,(n+1)s)}{k_P^\mathcal{M}(x,y,ns)}
\end{equation}
exists  for every $s>0$ (i.e., the ratio limit exists for every
``skeleton" sequence of the form $t_n=ns$, where $n=1,2,\dots$ and
$s>0$). Then
\begin{equation}\label{eqconjDw9}
\lim_{t\to\infty}\frac{k_P^\mathcal{M}(x,y,t+r)}{k_P^\mathcal{M}(x,y,t)}=1
\qquad \forall  r\in \mathbb{R}.
\end{equation}
\end{corollary}
\begin{proof}
 By Lemma \ref{assliminfsup}, the limit in
(\ref{eqliminf5}) equals $1$. By induction, $$ \lim_{n\to\infty}\!
\frac{k_P^\mathcal{M}(x,y,ns+r)}{k_P^\mathcal{M}(x,y,ns)}=1,$$
where $r\!=\!qs$, and $q\!\in\! \mathbb{Q}$, which (by Harnack)
implies that it holds  for $\forall r\!\in\! \mathbb{R}$. Hence,
\cite[Theorem 2]{Ki63} implies (\ref{eqconjDw9}).
 \end{proof}
\begin{rem}\label{remellHar}
If there exist $x_0,y_0\in \mathcal{M}$ and $0<s_0<1$ such that
\begin{equation}\label{eqconjDw17}
M(x_0,y_0,s_0):=\limsup_{t\to\infty}\frac{k_P^\mathcal{M}(x_0,y_0,t+s_0)}{k_P^\mathcal{M}(x_0,y_0,t)}<\infty,
\end{equation}
then by the parabolic Harnack inequality, for all $x,y,z,w\in
K\subset\subset \mathcal{M}$, $t>1$, we have the following Harnack
inequality of elliptic type:
\begin{equation*}\label{eqstronghara}
k_P^\mathcal{M}(z,w,t)\!\leq\!
C_1k_P^\mathcal{M}(x_0,y_0,t\!+\!\frac{s_0}{2})\! \leq\!
C_2k_P^\mathcal{M}(x_0,y_0,t\!-\!\frac{s_0}{2})\!\leq\!
C_3k_P^\mathcal{M}(x,y,t).
\end{equation*}
This estimate implies that for all $x,y\in \mathcal{M}$ and $r\in
\mathbb{R}$: \bean 0<m(x,y,r):=\liminf_{t\to\infty}
\frac{k_P^\mathcal{M}(x,y,t+r)}{k_P^\mathcal{M}(x_0,y_0,t)}
\leq\\\limsup_{t\to\infty}\frac{k_P^\mathcal{M}(x,y,t+r)}{k_P^\mathcal{M}(x_0,y_0,t)}
=M(x,y,r)<\infty. \eean
 \end{rem}
\begin{lem}\label{lem6}
 (a) The following assertions are equivalent:

(i) For each $x,y\in \mathcal{M}$ there exists a sequence $s_j\to
0$ of negative numbers such that for all $j\geq 1$, and $s=s_j$,
we have \begin{equation}\label{eqconjDw27}
\lim_{t\to\infty}\frac{k_P^\mathcal{M}(x,y,t+s)}{k_P^\mathcal{M}(x,y,t)}=1.
\end{equation}

(ii) The ratio limit in (\ref{eqconjDw27}) exists for any $x,y\in
\mathcal{M}$ and $s\in \mathbb{R}$.

(iii) Any limit function  $u(x,y,s)$ of the quotients
$\frac{k_P^\mathcal{M}(x,y,t_n+s)}{k_P^\mathcal{M}(x_0,x_0,t_n)}$
with $t_n\!\to\! \infty$  does not depend on $s$ and has the form
$u(x,y)$, where  $u(\cdot,y)\!\in\! \mathcal{C}_{P}(\mathcal{M})$
for every $y\!\in\! \mathcal{M}$ and $u(x, \cdot)\! \in\!
\mathcal{C}_{P^*}(\mathcal{M})$ for every $x\!\in\! \mathcal{M}$.

\vskip 2mm

(b) If one assumes further (\ref{eqoned}), then Conjecture
\ref{conjD} holds true.
\end{lem}
 \begin{proof} \hskip -2mm
(a) By Lemma~\ref{assliminfsup}, if the limit in
(\ref{eqconjDw27}) exists, then it is $1$.

(i) $\Rightarrow$ (ii). Fix $x_0,y_0\in \mathcal{M}$, and take
$s_0<0$ for which  the limit (\ref{eqconjDw27}) exists. It follows
that any limit solution $u(x,y,s)\in
\mathcal{H}_P(\mathcal{M}\times \mathbb{R}_-)$ of a sequence
$\frac{k_P^\mathcal{M}(x,y,t_n+s)}{k_P^\mathcal{M}(x_0,y_0,t_n)}$
with $t_n\to \infty$  satisfies
$u(x_0,y_0,s\!+\!s_0)\!=\!u(x_0,y_0,s)$ for all $s\!<\!0$. So,
$u(x_0,y_0,\cdot)$ is $s_0$-periodic. It follow from our
assumption and the continuity of $u$ that $u(x_0,y_0,\cdot)$ is
the constant function. Since this holds for all $x,y\!\in\!
\mathcal{M}$ and $u$, it follows that (\ref{eqconjDw27}) holds for
any $x,y\!\in\! \mathcal{M}$ and $s\!\in\! \mathbb{R}$.

(ii) $\Rightarrow$ (iii). Fix $y\in \mathcal{M}$. By Remark
\ref{6Rem21}, any limit function $u$ of the sequence
$\frac{k_P^\mathcal{M}(x,y,t_n+s)}{k_P^\mathcal{M}(x_0,x_0,t_n)}$
 with $t_n\to \infty$ belongs to
$\mathcal{H}_P(\mathcal{M}\times \mathbb{R}_-)$. Since
\begin{equation}\label{eqfrac6}
\frac{k_P^\mathcal{M}(x,y,t+s)}{k_P^\mathcal{M}(x_0,x_0,t)}=
\frac{k_P^\mathcal{M}(x,y,t)}{k_P^\mathcal{M}(x_0,x_0,t)}
\frac{k_P^\mathcal{M}(x,y,t+s)}{k_P^\mathcal{M}(x,y,t)}\,,
\end{equation}
(\ref{eqconjDw27}) implies that such a $u$ does not depend on $s$.
Therefore, $u=u(x,y)$, where $u(\cdot,y)\in
\mathcal{C}_{P}(\mathcal{M})$ and $u(x, \cdot) \in
\mathcal{C}_{P^*}(\mathcal{M})$.

(iii) $\Rightarrow$ (i). Write
\begin{equation}\label{eqquot67}
\frac{k_P^\mathcal{M}(x,y,t+s)}{k_P^\mathcal{M}(x,y,t)}=
\frac{k_P^\mathcal{M}(x,y,t+s)}{k_P^\mathcal{M}(x_0,x_0,t)}\,
\frac{k_P^\mathcal{M}(x_0,x_0,t)}{k_P^\mathcal{M}(x,y,t)}\,.
\end{equation}
Let $t_n\to \infty$ be a sequence such that the sequence
$\frac{k_P^\mathcal{M}(x,y,t_n+s)}{k_P^\mathcal{M}(x_0,x_0,t_n)}$
converges to a solution in $\mathcal{H}_P(\mathcal{M}\times
\mathbb{R}_-)$. By our assumption, we have
$$\lim_{n\to\infty}\frac{k_P^\mathcal{M}(x,y,t_n+s)}{k_P^\mathcal{M}(x_0,x_0,t_n)}
=\lim_{n\to\infty}\frac{k_P^\mathcal{M}(x,y,t_n)}{k_P^\mathcal{M}(x_0,x_0,t_n)}=u(x,y)>0,$$
which together with (\ref{eqquot67}) imply (\ref{eqconjDw27}) for
all $s\in \mathbb{R}$.

\vspace{1 mm}

\noindent (b) The uniqueness and (iii) imply that
$\frac{k_P^\mathcal{M}(x,y,t+s)}{k_P^\mathcal{M}(x_0,x_0,t)}\rightarrow\!
\frac{u(x)u^*(y)}{u(x_0)u^*(x_0)}$, where $u\!\in\!
\mathcal{C}_{P}(\mathcal{M})$ and $u^*\!\in\!
\mathcal{C}_{P^*}(\mathcal{M})$, and Conjecture \ref{conjD} holds.
 \end{proof}
 \begin{rem}\label{rem15}
Let $\mathcal{M}\subsetneqq \mathbb{R}^d$ be a smooth domain and
$P$ and $P^*$ be (up to the boundary) smooth operators. Denote by
$\mathcal{C}^0_{P}(\mathcal{M})$ the cone of all functions in
$\mathcal{C}_{P}(\mathcal{M})$ which vanish on $\partial
\mathcal{M}$. Suppose that one of the conditions (i)--(iii) of
Lemma \ref{lem6} is satisfied. Clearly, for any fixed $y$ any
limit function $u(\cdot,y)$ of Lemma \ref{lem6} belongs to the
Martin boundary `at infinity' which in this case is
$\mathcal{C}^0_{P}(\mathcal{M})$. Therefore, Conjecture
\ref{conjD} holds true if the Martin boundaries `at infinity' of
$P$ and $P^*$ are one-dimensional. As a simple example, take
$P=-\Delta$ and $\mathcal{M}=\mathbb{R}^d_+$.
 \end{rem}
\begin{lem}\label{lem63}
Suppose that $P$ is null-critical, and for each $x,y\in
\mathcal{M}$ there exists a sequence  $\{s_j\}$  of negative
numbers such that $s_j \to 0$, and
\begin{equation}\label{eqconjDw71}
\liminf
_{t\to\infty}\frac{k_P^\mathcal{M}(x,y,t+s)}{k_P^\mathcal{M}(x,y,t)}\geq
1
\end{equation}
for $s=s_j$, $j=1,2,\ldots\,$. Then Conjecture \ref{conjD} holds
true.
\end{lem}
 \begin{proof} Let $u(x,y,s)$ be a limit function of a sequence
$\frac{k_P^\mathcal{M}(x,y,t_n+s)}{k_P^\mathcal{M}(x_0,x_0,t_n)}$
with $t_n\to \infty$ and $s<0$.  By our assumption,
 $u(x,y,s+s_j)\geq u(x,y,s)$, and therefore,
$u_s(x,y,s)\leq 0$ for all $s<0$. Thus, $u(\cdot,y,s)$ (respect.,
$u(x,\cdot,s)$) is a positive supersolution of the equation $Pu=0$
(respect., $P^*u=0$) in $\mathcal{M}$. Since $P$ is critical, it
follows that $u(\cdot,y,s)\in \mathcal{C}_{P}(\mathcal{M})$
(respect., $u(x,\cdot,s)\in \mathcal{C}_{P^*}(\mathcal{M})$), and
hence $u_s(x,y,s)=0$. By the uniqueness, $u$ equals to
$\frac{\varphi(x)\varphi^*(y)}{\varphi(x_0)\varphi^*(x_0)}$, and
Conjecture \ref{conjD} holds true.
 \end{proof}
\begin{rem}\label{questTaub} Suppose that $P$ is null-critical,
and fix $x_0\neq y_0$. Then using Theorem~\ref{mainthm} and
\cite[Theorem~2.1]{P90} we have for $x\neq y$:
$$\mbox{ \em{(i)} }\lim_{t\to \infty} k_P^\mathcal{M}(x,y,t)=\lim_{t\to \infty} k_P^\mathcal{M}(x_0,y_0,t)=0,$$
  $$\mbox{ \em{(ii)} }\int_0^\infty k_P^\mathcal{M}(x,y,\tau)\,\mathrm{d}\tau=
\int_0^\infty
k_P^\mathcal{M}(x_0,y_0,\tau)\,\mathrm{d}\tau=\infty,$$
  $$ \mbox{ \em{(iii)} }\lim_{\lambda\nearrow 0}\!\frac{\int_0^\infty\!
\mathrm{e}^{\lambda\tau}k_P^\mathcal{M}(x,y,\tau)\mathrm{d}\tau}
{\int_0^\infty\! \mathrm{e}^{\lambda\tau}k_P^\mathcal{M}(x_0,y_0,
\tau) \mathrm{d}\tau}=\lim_{\lambda\nearrow
0}\frac{\Green{\mathcal{M}}{P-\lambda}{x}{y}}
{\Green{\mathcal{M}}{P-\lambda}{x_0}{y_0}} =
\frac{\varphi(x)\varphi^*(y)}{\varphi(x_0)\varphi^*(y_0)}.$$
Therefore, Conjecture \ref{conjD} would follow from a strong ratio
Tauberian theorem if additional Tauberian conditions are satisfied
(see, \cite{BGT,ST}).
\end{rem}
\section{The parabolic Martin boundary} The large time behavior
of quotients of the heat kernel is obviously closely related to
the parabolic Martin boundary (for the parabolic Martin boundary
theory see \cite{Doo}). Theorem \ref{thmconjD} relates Conjecture
\ref{conjD} and the parabolic Martin compactification of
$\mathcal{H}_P(\mathcal{M}\times \mathbb{R}_-)$.

\begin{lem}\label{thmconjD7}
Fix $y\in \mathcal{M}$. The following assertions are equivalent:

(i) For each $x\in \mathcal{M}$ there exists a sequence $s_j\to 0$
of negative numbers such that
\begin{equation}\label{eqconjDw77}
\lim_{t\to\infty}\frac{k_P^\mathcal{M}(x,y,t+s)}{k_P^\mathcal{M}(x,y,t)}
\end{equation}
exists for $s=s_j$, $j=1,2,\ldots\,$.

(ii) Any parabolic Martin function in
$\mathcal{H}_P(\mathcal{M}\times \mathbb{R}_-)$ corresponding to a
Martin sequence $\{(y,-t_n)\}_{n\!=\!1}^\infty$, where $t_n\!\to\!
\infty$, is time independent.
\end{lem}
  \begin{proof}   Let $K_P^\mathcal{M}(x,y,s)=\lim_{n\to\infty}\frac{k_P^\mathcal{M}(x,y,t_n+s)}
{k_P^\mathcal{M}(x_0,y,t_n)}$ be such a Martin function. The lemma
follows from the identity
$$\frac{k_P^\mathcal{M}(x,y,t_n+s)}{k_P^\mathcal{M}(x_0,y,t_n)}=
\frac{k_P^\mathcal{M}(x,y,t_n+s)}{k_P^\mathcal{M}(x,y,t_n)}\,
\frac{k_P^\mathcal{M}(x,y,t_n)}{k_P^\mathcal{M}(x_0,y,t_n)}\,,$$
and Lemma \ref{lem6}.
 \end{proof}
\begin{thm}\label{thmconjD}
Assume that (\ref{eqconjDw17}) holds true for some $x_0,y_0\in
\mathcal{M}$, and $s_0>0$. Then the following assertions are
equivalent:

(i) Conjecture \ref{conjD} holds true for a fixed $x_0 \in
\mathcal{M}$.

  (ii)
\begin{equation}\label{eqconjD2}
\lim_{t\to\infty}\frac{k_P^\mathcal{M}(x,y,t)}{k_P^\mathcal{M}(x_1,y_1,t)}
\end{equation}
exists, and the limit is positive for every $x,y,x_1,y_1\in
\mathcal{M}$.

  (iii)
\begin{equation}\label{eqconjD1}
\lim_{t\to\infty}\frac{k_P^\mathcal{M}(x,y,t)}{k_P^\mathcal{M}(y,y,t)}
\, , \quad \mbox{and} \quad
\lim_{t\to\infty}\frac{k_P^\mathcal{M}(x,y,t)}{k_P^\mathcal{M}(x,x,t)}
\end{equation}
exist, and these ratio limits are positive for every $x,y\in
\mathcal{M}$.

(iv) For any $y\in \mathcal{M}$ there is a unique nonzero
parabolic Martin boundary point $\bar y$ for the equation  $Lu=0$
in $\mathcal{M}\times \mathbb{R}$ which corresponds to any
sequence of the form $\{(y,-t_n)\}_{n=1}^\infty$ such that $t_n\to
\infty$, and for any $x\in \mathcal{M}$ there is a unique nonzero
parabolic Martin boundary point $\bar x$ for the equation
$u_t+P^*u=0$ in $\mathcal{M}\times \mathbb{R}$ which corresponds
to any sequence of the form $\{(x,-t_n)\}_{n=1}^\infty$ such that
$t_n\to \infty$.

Moreover, if Conjecture \ref{conjD} holds true, then for any fixed
$y\in \mathcal{M}$ (respect.,  $x \in \mathcal{M}$), the limit
function $a(\cdot,y)$ (respect., $a(x,\cdot)$) is a positive
solution of the equation $Pu = 0$ (respect., $P^*u = 0$).
Furthermore, the Martin functions of part (iv) are time
independent, and (\ref{eqconjDw27}) holds for all $x,y \in
\mathcal{M}$ and $s \in \mathbb{R}$.
\end{thm}
  \begin{proof}  (i) $\Rightarrow$ (ii) follows from the identity
$$\frac{k_P^\mathcal{M}(x,y,t)}{k_P^\mathcal{M}(x_1,y_1,t)}=
\frac{k_P^\mathcal{M}(x,y,t)}{k_P^\mathcal{M}(x_0,x_0,t)}\cdot
\left(\frac{k_P^\mathcal{M}(x_1,y_1,t)}{k_P^\mathcal{M}(x_0,x_0,t)}\right)^{-1}.$$

(ii) $\Rightarrow$ (iii). Take $x_1=y_1=y$ and $x_1=y_1=x$,
respectively.

(iii) $\Rightarrow$ (iv). It is well known that the Martin
compactification does not depend on the fixed reference point
$x_0$. So, fix $y\in \mathcal{M}$ and take it also as a reference
point. Let $\{-t_n\}$ be a sequence such that $t_n\to \infty$ and
such that the Martin sequence
$\frac{k_P^\mathcal{M}(x,y,t+t_n)}{k_P^\mathcal{M}(y,y,t_n)}$
converges to a Martin function $K_P^\mathcal{M}(x,\bar{y},t)$. By
our assumption, for any $t$ we have
$$\lim_{n\to\infty}\frac{k_P^\mathcal{M}(x,y,t+t_n)}{k_P^\mathcal{M}(y,y,t+t_n)}=
\lim_{\tau\to\infty}\frac{k_P^\mathcal{M}(x,y,\tau)}{k_P^\mathcal{M}(y,y,\tau)}=
b(x)>0,$$ where $b$ does not depend on the sequence $\{-t_n\}$. On
the other hand,
$$\lim_{n\to\infty}\frac{k_P^\mathcal{M}(y,y,t+t_n)}{k_P^\mathcal{M}(y,y,t_n)}=
K_P^\mathcal{M}(y,\bar{y},t)=f(t).$$
 Since
$$\frac{k_P^\mathcal{M}(x,y,t+t_n)}{k_P^\mathcal{M}(y,y,t_n)}=
\frac{k_P^\mathcal{M}(x,y,t+t_n)}{k_P^\mathcal{M}(y,y,t+t_n)}\cdot
\frac{k_P^\mathcal{M}(y,y,t+t_n)}{k_P^\mathcal{M}(y,y,t_n)},$$ we
have
$$K_P^\mathcal{M}(x,\bar{y},t)=b(x)f(t).$$
 By separation of variables, there exists a constant
$\lambda$ such that
$$Pb-\lambda b=0 \quad \mbox{ on } \mathcal{M},
\qquad f'+\lambda f=0 \quad \mbox{ on } \mathbb{R},\;\; f(0)=1.
$$
 Since $b$ does not depend on the sequence
$\{-t_n\}$, it follows in particular,  that $\lambda$ does not
depend on this sequence. Thus, $ \lim_{\tau\to \infty}
\frac{k_P^\mathcal{M}(y,y,t+\tau)}{k_P^\mathcal{M}(y,y,\tau)}=f(t)=\mathrm{e}^{-\lambda
t}$. Lemma~\ref{assliminfsup} implies that $\lambda=0$.  It
follows that $b$ is a positive solution of the equation $Pu=0$,
and
\begin{equation}\label{6eq1}
K_P^\mathcal{M}(x,\bar{y},t)=\lim_{\tau\to -
\infty}\frac{k_P^\mathcal{M}(x,y,t-\tau)}{k_P^\mathcal{M}(y,y,-\tau)}=
b(x).
\end{equation}
The dual assertion can be proved similarly.

(iv) $\Rightarrow$ (i). Let $K_P^\mathcal{M}(x,\bar{y},t)$ be a
Martin function, and $s_0>0$ such that
$K_P^\mathcal{M}(x_0,\bar{y},s_0/2)>0$. Consequently,
$K_P^\mathcal{M}(x,\bar{y},s)>0$  for $s\geq s_0$. Using the
substitution $\tau=s+s_0$ we obtain \bean \!\!\!\!\!\lim_{\tau\to
\infty}\!
\frac{k_P^\mathcal{M}(x,y,\tau)}{k_P^\mathcal{M}(x_0,x_0,\tau)}\!=\!
\lim_{s\to \infty}\!\!\left \{\!\!
\frac{k_P^\mathcal{M}(x,y,s+s_0)}{k_P^\mathcal{M}(y,y,s)}\right.\times\\[3mm]
\left.\frac{k_P^\mathcal{M}(y,y,s)}{k_P^\mathcal{M}(x_0,y,s\!+\!2s_0)}
\!\frac{k_P^\mathcal{M}(x_0,y,s\!+\!2s_0)}{k_P^\mathcal{M}(x_0,x_0,s\!+\!s_0)}\!\!\right
\}\!=
\!\!\frac{K_P^\mathcal{M}(x,\bar{y},s_0)K_{P^*}^\mathcal{M}(\overline{x_0},y,s_0)}
{K_P^\mathcal{M}(x_0,\bar{y},2s_0)}.
 \eean
The last assertion of the lemma follows from  (\ref{6eq1}) and
Lemma \ref{lem6}.
 \end{proof}
 \section{Minimal positive solutions}
In this section we discuss the relation between Conjecture
\ref{conjD} and the parabolic and elliptic {\em minimal} Martin
boundaries.
\begin{rem}\label{remmin8}
By the parabolic Harnack inequality for $P^*$, we have for each
$0<\varepsilon<1$
\begin{equation}\label{eqsrhi}
k_P^\mathcal{M}(x,y_0,t-\varepsilon)\leq C(y_0,\varepsilon)
k_P^\mathcal{M}(x,y_0,t) \qquad \forall x\in \mathcal{M},\; t>1.
\end{equation}
 Therefore, if
$\{(y_0,t_n)\}$ is a nontrivial minimal Martin sequence with
$t_n\to -\infty$, then one infers as in \cite{KT} that the
corresponding minimal parabolic function in
$\mathcal{H}_P(\mathcal{M}\times \mathbb{R}_-)$ is of the form
$u(x,t)= \mathrm{e}^{-\lambda t}u_\lambda(x,y_0)$ with
$\lambda\leq 0$ and $u_\lambda\in
\mathrm{exr\,}\mathcal{C}_{P-\lambda}(\mathcal{M})$, where
$\mathrm{exr\,}\mathcal{C}$ is the set of extreme rays of a cone
$\mathcal{C}$. If further, for some $x_0\in \mathcal{M}$ and $s<0$
one has
\begin{equation}\label{eqlimitexist}
\liminf_{t\to\infty}\frac{k_P^\mathcal{M}(x_0,y_0,t+s)}{k_P^\mathcal{M}(x_0,y_0,t)}\geq
1,
\end{equation}
then $\lambda=0$, and consequently, $u$ is also a minimal solution
in $\mathcal{C}_{P}(\mathcal{M})$. Recall that in the selfadjoint
case, the ratio {\em limit} in (\ref{eqlimitexist}) equals $1$.
\end{rem}
\begin{lem}\label{cormin}
Suppose that the ratio limit in (\ref{eqconjDw27}) exists for all
$x,y\in \mathcal{M}$ and $s\in \mathbb{R}$. Let
$a(x,y):=\lim_{n\to
\infty}\frac{k_P^\mathcal{M}(x,y,t_n+s)}{k_P^\mathcal{M}(x_0,x_0,t_n)}$,
where $t_n\!\to\! \infty$. If for some $y_0\in \mathcal{M}$ the
function $u(x):=a(x,y_0)$ is minimal in
$\mathcal{C}_{P}(\mathcal{M})$, then $a(x,y)=u(x)v(y)$, where
$v\in\mathcal{C}_{P^*}(\mathcal{M})$.
\end{lem}
\begin{proof}
Fix $y\in \mathcal{M}$ and $\varepsilon>0$. By the parabolic
Harnack inequality for $P^*$ and Lemma \ref{lem6}, we have
\begin{equation}\label{eqcormin}
\frac{k_P^\mathcal{M}(x,y,t-\varepsilon)}{k_P^\mathcal{M}(x_0,x_0,t)}
\leq
C(y,\varepsilon)\frac{k_P^\mathcal{M}(x,y_0,t)}{k_P^\mathcal{M}(x_0,x_0,t)}\qquad
\forall x\in \mathcal{M} .
\end{equation}
Therefore, $a(x,y)\leq C(y)u(x)$ which implies the claim.
\end{proof}
The following examples demonstrate that if Conjecture \ref{conjD}
holds true while (\ref{eqoned}) does not hold, then the limit
function $a(\cdot,y)$ is typically a non-minimal solution in
$\mathcal{C}_P(\mathcal{M})$.
\begin{example}\label{examinim}
Consider a (regular) Benedicks domain $\mathcal{M}\!\subseteq\!
\mathbb{R}^d$ such that the cone of positive harmonic functions
which vanish on $\partial \mathcal{M}$ is of dimension two. By
\cite{CMS1}, Conjecture \ref{conjD} holds true in this case, the
limit function is not a product of two (separated) harmonic
functions,  and therefore, $a(\cdot,y)$ is not minimal in
$\mathcal{C}_{-\Delta}(\mathcal{M})$ for any $y\in \mathcal{M}$.
 \end{example}
\begin{example}\label{examinim1}
Consider a radially symmetric Schr\"{o}dinger operator
$H\!:=\!-\Delta\!+\!V(|x|)$ on $\mathbb{R}^d$ with a bounded
potential. Suppose that $\lambda_0=0$, and that the Martin
boundary of $H$ on $\mathbb{R}^d$ is homeomorphic to $S^{d-1}$
(see \cite{M86}). Clearly, any Martin function corresponding to
$\{(y_0,t_n)\}$ with $x_0\!=\!y_0\!=\!0$ is radially symmetric. It
follows that Davies' conjecture holds true for $x_0\!=\!y\!=\!0$,
and the limit function is the normalized positive radial solution
in $\mathcal{C}_{H}(\mathbb{R}^d)$. This solution is not minimal
in $\mathcal{C}_{H}(\mathbb{R}^d)$. Thus, any limit function
$u(\cdot,y)$ is not minimal in $\mathcal{C}_{H}(\mathbb{R}^d)$.
 \end{example}
We conclude this section with some related problems. The following
conjecture was posed by the author in \cite[Conjecture
3.6]{Pheat}.
\begin{conjecture}\label{conjP}
Suppose that $P$ is critical operator in $\mathcal{M}$, then the
ground state $\varphi$ is a minimal positive solution in the cone
$\mathcal{H}_P(\mathcal{M}\times \mathbb{R})$.
\end{conjecture}
Note that if (\ref{eqconjDw71}) holds true, then by Theorem
\ref{thmconjD}, the ground state is a Martin function in
$\mathcal{H}_P(\mathcal{M}\times \mathbb{R})$.

\begin{example}\label{question4} Consider again the example in \cite[Section
4]{PSt}. In that example $-\Delta$ is subcritical in
$\mathcal{M}$, $\lambda_0=0$, and (\ref{eqoned}) and
Conjecture~\ref{eqoned} hold true. Hence, $\mathbf{1}$ is a Martin
function in $\mathcal{H}_{-\Delta}(\mathcal{M}\times \mathbb{R})$.
On the other hand, $\mathbf{1}\in
\mathrm{exr\,}\mathcal{C}_{-\Delta}(\mathcal{M})$ but
$\mathbf{1}\not\in
\mathrm{exr\,}\mathcal{H}_{-\Delta}(\mathcal{M}\times
\mathbb{R})$. So, Conjecture~\ref{conjP} cannot be extended to the
subcritical ``Liouvillian" case (see also \cite{BS}).
\end{example}
Thus, it would be interesting to study the following problem which
was raised by Burdzy and Salisbury \cite{BS} for $P=-\Delta$ and
$\mathcal{M}\subset \mathbb{R}^d$.
\begin{question}\label{questBS}
Assume that $\lambda_0=0$. Determine which minimal positive
solutions in $\mathcal{C}_{P}(\mathcal{M})$ are minimal in
$\mathcal{H}_P(\mathcal{M}\times \mathbb{R_-})$.
\end{question}
\section{Uniform Harnack
inequality and Davies' conjecture}\label{securhi}
 In this section we discuss the relationship between the parabolic
Martin boundary of $\mathcal{H}_P(\mathcal{M}\times
\mathbb{R}_-)$, the elliptic Martin boundaries of
$\mathcal{C}_{P-\lambda}(\mathcal{M})$, $\lambda\leq\lambda_0=0$,
and Conjecture \ref{conjD} under a certain assumption.
\begin{defn} \label{DUHI}
We say that the {\em uniform restricted parabolic Harnack
inequality} (in short, (URHI)) holds in
$\mathcal{H}_P(\mathcal{M}\times \mathbb{R}_-)$ if
for any $\varepsilon>0$ there exists a positive constant
$C=C(\varepsilon)>0$ such that
 \begin{equation} \label{UHI}
 u(x,t-\varepsilon) \!\leq\!
Cu(x,t)\quad  \forall(x,t) \!\in\!
\mathcal{M}\!\times\!\mathbb{R}_-\mbox{ and } \forall u \!\in\!
\mathcal{H}_P(\mathcal{M}\!\times\! \mathbb{R}_-).
 \end{equation}
 \end{defn}
 It is well known that (URHI) holds true if and only
if the {\em separation principle} (SP) holds true, that is, $u\neq
0$ is in $\mathrm{exr\,}\mathcal{H}_P(\mathcal{M}\times
\mathbb{R}_-)$ if and only if $u$ is of the form
$\mathrm{e}^{-\lambda t}v_\lambda(x)$, where $v_\lambda\in
\mathrm{exr\,}\mathcal{C}_{P-\lambda}(\mathcal{M})$ \cite{KT,M}.
In particular, the answer to Question \ref{questBS} is simple if
(URHI) holds.
\begin{lem}\label{lemurhi1}
(i) Suppose that (URHI) holds true, then for any $s<0$
$$\ell_+:=\limsup_{t\to\infty}
\frac{k_P^\mathcal{M}(x,y,t+s)}{k_P^\mathcal{M}(x,y,t)}\leq 1
\qquad\mbox{(Lin's condition)}.$$

 (ii) Assume further
that for some  $x_0, y_0 \in \mathcal{M}$ and $s_0<0$
$$\ell_-:=\liminf_{t\to\infty}
\frac{k_P^\mathcal{M}(x_0,y_0,t+s_0)}{k_P^\mathcal{M}(x_0,y_0,t)}
\geq 1 ,$$  then any limit function $u(x,y,s)$ of $
\frac{k_P^\mathcal{M}(x,y,t_n+s)}{k_P^\mathcal{M}(x_0,y_0,t_n)}$
with $t_n\to \infty$ does not depend on $s$, and has the form
$u(x,y)$, where  $u(\cdot,y)\in \mathcal{C}_{P}(\mathcal{M})$  for
every $y\in \mathcal{M}$ and $u(x, \cdot) \in
\mathcal{C}_{P^*}(\mathcal{M})$ for every $x\in \mathcal{M}$.

(iii) If one assumes further (\ref{eqoned}), then Conjecture
\ref{conjD} holds true.
\end{lem}
\begin{proof} (i) By (URHI), if
 $u\in \mathrm{exr\,}\mathcal{H}_P(\mathcal{M}\times \mathbb{R}_-)$, then
$u(x,t)=\mathrm{e}^{-\lambda t}u_\lambda(x)$, where $\lambda\leq
0$. Consequently, for every $u\in \mathcal{H}_P(\mathcal{M}\times
\mathbb{R}_-)$
\begin{equation}\label{equrhi}
u(x,t+s)\leq u(x,t) \qquad \forall(x,t)\in \mathcal{M}\!\times\!
\mathbb{R}_-,\mbox{ and } \forall s<0,
\end{equation}
and equality holds for some $\!s<0\!$ and $(x,t)\!\in\!
\mathcal{M}\!\times\! \mathbb{R}_-$ if and only if $u\!\in\!
\mathcal{C}_{P}(\mathcal{M})$. Clearly, (\ref{equrhi}) implies
that  $$\ell_+:=\limsup_{t\to\infty}
\frac{k_P^\mathcal{M}(x,y,t+s)}{k_P^\mathcal{M}(x,y,t)} \leq 1
\quad\forall x, y \in \mathcal{M} \mbox{ and } s<0,$$ which
together with Lemma \ref{assliminfsup} imply $\ell_+=1$.

(ii) At the point $(x_0,y_0,s_0)$ we have $\ell_-=\ell_+=1$,
therefore,
\begin{equation}\label{eqskeleton}
\lim_{t\to\infty}
\frac{k_P^\mathcal{M}(x_0,y_0,t+s_0)}{k_P^\mathcal{M}(x_0,y_0,t)}=1.
\end{equation}
Consequently, for any sequence $t_{k}\to\infty$ satisfying
$$\lim_{k\to\infty}
\frac{k_P^\mathcal{M}(x,y_0,t_{k}+\tau)}{k_P^\mathcal{M}(x_0,y_0,t_{k})}=u(x,\tau)\quad
\forall (x,\tau)\in \mathcal{M} \times \mathbb{R}_- ,$$ we have
$u(x_0,s_0)=u(x_0,2s_0)=1$, and therefore, $u\in
\mathcal{C}_{P}(\mathcal{M})$. The other assertions of the lemma
follow from Lemma \ref{lem6}. \end{proof}
\begin{rem}\label{remurhi3}
From the proof of Lemma \ref{lemurhi1} it follows that if (URHI)
holds true, then a sequence $t_n\to \infty$ satisfies
$$\lim_{n\to\infty}
\frac{k_P^\mathcal{M}(x_0,y_0,t_n+s_0)}{k_P^\mathcal{M}(x_0,y_0,t_n)}
=1,$$ for some $x_0, y_0 \in \mathcal{M}$ and $s_0\neq 0$ if and
only if
$$\lim_{n\to\infty}
\frac{k_P^\mathcal{M}(x,y,t_n+s)}{k_P^\mathcal{M}(x,y,t_n)}
=1\qquad \forall x, y \in \mathcal{M} \mbox{ and } s\in
\mathbb{R}.$$

 \end{rem}
\begin{corollary}\label{cormartin2}
Suppose that (URHI) holds true, then there exists a sequence
$t_n\to \infty$ such that
$\lim_{n\to\infty}\frac{k_P^\mathcal{M}(x,y,t_n)}{k_P^\mathcal{M}(x_0,x_0,t_n)}=a(x,y)
$
exists and is positive for all $x,y\in \mathcal{M}$. Moreover,
$a(\cdot,y)\in\mathcal{C}_{P}(\mathcal{M})$, and $a(\cdot,y)$ is a
parabolic Martin function for all $y\in \mathcal{M}$. For each
$x\in \mathcal{M}$ the function $a(x,\cdot)$ satisfies similar
properties with respect to $P^*$.
\end{corollary}
\begin{proof} Take $s_0\neq 0$ and $\{t_n\}$ such that
$\lim_{n\to\infty}
\frac{k_P^\mathcal{M}(x_0,y_0,t_n+s_0)}{k_P^\mathcal{M}(x_0,y_0,t_n)}
=1,$ and use Remark \ref{remurhi3} and a standard diagonalization
argument. \end{proof}

\end{document}